\documentclass[11pt, oneside]{article}   	
\usepackage{geometry}                		
\usepackage{graphicx}
\usepackage{amssymb}
\usepackage{amsthm}

\newtheorem{theorem}{Theorem}

\newtheorem{cor}{Corollary}

\title{Betti numbers of fat forests and their Alexander dual}
\author{Ralf Fr\"oberg\footnote{Mathematics Department Stockholm University Sweden email:frobergralf@gmail.com, ORCID 0000-0002-7294-2856}}
\date{}							% Activate to display a given date or no date

\begin{document}
\maketitle
%\section{}
%\subsection{}
\begin{abstract}
Let $k$ be a field and $R=k[x_1,\ldots,x_n]/I=S/I$ a graded ring. Then $R$ has a $t$-linear
resolution if $I$ is generated by homogeneous elements of degree $t$, and all higher syzygies
are linear. Thus $R$ has a $t$-linear resolution if ${\rm Tor}^S_{i,j}(S/I,k)=0$ if $j\ne i+t-1$.

For a graph $G$ on $\{1,\ldots,n\}$, the edge algebra is $k[x_1,\ldots,x_n]/I$, where $I$ is generated by those
$x_ix_j$ for which $\{ i,j\}$ is an edge in $G$. We want to determine the Betti numbers of edge rings with
2-linear resolution. But we want to do that by looking at the edge ring as a Stanley-Reisner ring.

For a simplicial complex $\Delta$ on $[{\bf n}]=\{1,\ldots,n\}$ and a field $k$, the Stanley-Reisner
 ring $k[\Delta]$ is $k[x_1,\ldots,x_n]/I$, where $I$ is generated by the squarefree monomials
 $x_{i_1}\cdots x_{i_k}$ for which $\{ i_1,\ldots,i_k\}$ does not belong to $\Delta$. 

Which Stanley-Reisner rings that are edge rings with 2-linear resolution are known. Their
associated complexes has had different names in the literature. We call them fat forests here.
We determine the Betti numbers of many fat forests and compare our result with what is known. We also consider Betti numbers
of Alexander duals of fat forests.

Keywords: Stanley-Reisner ring, edge ring, Betti numbers, Hilbert series
\end{abstract}

\section{Background}

The simplicial complexes we will consider have had different names.
They are called Generalized forests (\cite{Fr}), Quasiforests (\cite{Zh}), or Fat forests (\cite{Au}). 
We will call them fat forests.
They are recursively defined as follows.
A $d$-simplex $F_1$ of dimension $\ge0$ (i.e. with $d+1$ vertices) 
is a fat forest. If $F_i$, $i=1,\ldots,k$, are simplices and
$G_{k-1}=F_1\cup\cdots\cup F_{k-1}$ is a fat forest, then $G_{k-1}\cup F_k$ is a fat forest if 
$H=G_{k-1}\cap F_k$ is a simplex, $\dim H\ge -1$. (If $\dim H=-1$, 
then $G_{k-1}$ and $F_k$ are disjoint.) 
 A fat forest a called a fat tree if it is connected, so if $\dim G_{k-1}\cap F_k\ge0$ for all $k$,
 but here it is not is necessary to treat fat trees separately.
 
 \bigskip
 Let $G$ be a graph on $[{\bf n}]$ and let $S=k[x_1,\ldots,x_n]$, where $k$ is a field. 
The edge ring of $G$ is $S/I$,where $I$ is generated by all $x_ix_j$ for which $\{ i,j\}$ is an edge 
of $G$. Let $G^c$ be the complement graph to $G$, i.e. the graph on $[{\bf n}]$ with edges
$\{ k,l\}$ for which $\{ k,l\}$ is not an edge of $G$. It is shown in \cite{Fr} that the edge ring of $G$
has a 2-linear resolution, meaning ${\rm Tor}^S_{i,j}(S/I,k)=0$ if $j\ne i+1$, if and only if 
the complementary graph $G^c$ (having those edges which are not edges in $G$) is 
chordal. A graph is chordal if every $k$-cycle, $k\ge4$, has a chord. This theorem has been reproved
in different ways in \cite{E-G-H-P}, \cite{He-Hi}, \cite{En-St}, and \cite{vT}.
 An edge ring $k[x_1,\ldots,x_n]/I$ can also be seen as a Stanley-Reisner ring $k[\Delta]$. Then
$I$ is generated by those squarefree monomials that do not belong to a simplicial complex $\Delta$.
The 1-skeleton of a
simplicial complex $\Delta$ consists of all faces of $\Delta$ that have dimension $\le1$.
Dirac has shown in \cite[Theorem 1 and 2]{Di}, that a graph $G$ is chordal if and only if it is the 1-skeleton of a fat forest. There is an algebraic proof of Dirac's theorems in \cite{Zh1}. 
Dirac's theorem gives
easily, see \cite{Fr}, that the Stanley-Reisner ring $k[\Delta]$ has a 2-linear resolution if and only
if $\Delta$ is a fat forest. 

\bigskip
Edge rings of Ferrer's graphs have 2-linear resolution. Their resolution has 
been determined in \cite{Co-Na}. Other classes of monomial rings with 2-linear resolution are treated in e.g. \cite{Co-Na1}, \cite{Ho},
\cite{Ar-He-Hi}, and \cite{Ja}. Finally  minimal resolution of all edge rings with 2-linear 
resolution was determined
in \cite{Ch}. This gives a (rather complicated)
way to determine the Betti numbers. In this note we will show that if one is only interested
in the Betti numbers of edge rings with 2-linear resolution, there is in many cases a very easy way to determine them using the description
of the rings as Stanley-Reisner rings.

\bigskip
Monomial rings with 2-linear resolution are also treated in e.g. \cite{Ng},
\cite{Wo}, and \cite{He-Hi-Zh}.

\section{Hilbert series and Betti numbers of fat forests}
If $S/I$ has a 2-linear resolution it looks like this:

$$0\leftarrow S/I\leftarrow S\leftarrow S[-2]^{b_1}\leftarrow S[-3]^{b_2}\leftarrow\cdots\leftarrow S[-p-1]^{b_p}\leftarrow 0$$

where $S[-i]$ means that we have shifted degrees of $S$ $i$ steps.

\bigskip
Using that the alternating sum of the $k$-dimensions in each degree is 0, we get that the
Hilbert series of $k[\Delta]$ with 2-linear resolution equals 
$\frac{1-\beta_{1,2}t^2+\beta_{2,3}t^3-\cdots(-1)^{p}\beta_{p,p+1}t^{p+1}}{(1-t)^n}$, where $\beta_{i,j}$ are 
the graded Betti numbers $\dim_k{\rm Tor}^S_{i,j}(k[\Delta],k)$, and $n$ is the
number of vertices in $\Delta$.
We are primarily interested in the Betti numbers $\beta_{i,j}=\dim_k{\rm Tor}^S_{i,j}(S/I,k)$ of
Stanley-Reisner rings of fat trees, but the Hilbert series
contains the same information as the set of Betti numbers.

\bigskip
\noindent
{\em Example} Let $\Delta$ be the simplicial complex with facets (maximal faces) $\{1,2\}$,
$\{2,3,4\}$, and $\{5\}$. This is a fat forest which can be built in the following way. Start with the simplex $\{2,3,4\}$ 
which has Hilbert series $\frac{1}{(1-t)^3}$. Then adjoin $\{1,2\}$ in $\{2\}$. We add the Hilbert series $\frac{1}{(1-t)^2}$
for the face $\{1,2\}$ and subtract $\frac{1}{1-t}$ for the face $\{2\}$ which has been counted twice. Finally adjoin $\{5\}$
with Hilbert series 1. Thus the Hilbert series of $\Delta$ is $\frac{1}{(1-t)^3}+\frac{1}{(1-t)^2}-\frac{1}{1-t}+1=
\frac{1+2t-3t^2+t^3}{(1-t)^3}=\frac{1-6t^2+9t^3-5t^4+t^5}{(1-t)^5}$, and 
the Betti numbers of $k[\Delta]$
are $b_{1,2}=6,b_{2,3}=9,b_{3,4}=5$, and $b_{4,5}=1$. This ring can also be seen as the edge ring of a graph with
edges \{1,3\}, \{1,4\}, \{1,5\}, \{2,5\}, \{3,5\}, and \{4,5\}.
 
 \begin{theorem} Let $F=F_1\cup\cdots\cup F_k$ be a fat tree with $F_i$ a simplex of dimension $d_i$
 and $(F_1\cup\cdots\cup F_{j-1})\cap F_j$ a simplex of dimension $r_j$. 
 Then the Hilbert series of $k[F]$
 is $\sum_{i=1}^k\frac{1}{(1-t)^{d_i+1}}-\sum_{i=2}^{k}\frac{1}{(1-t)^{r_i+1}}$.
 The projective dimension is $\sum_{i=1}^kd_i-\sum_{i=2}^kr_i+1-\min\{ r_i\}-2$. 
 The depth of $k[F]$ is $\min\{ r_i\}+2$, and $F$ is CM (Cohen-Macaulay) 
 if and only if there is a $d$ such that
 $d_i=d$ for all $i$ and $r_i=d-1$ for all $i$.
  \end{theorem}
 
 \noindent
 {\em Proof} The definition of fat forests directly gives the Hilbert series.
The number of vertices of $F$ is $\sum_{i=1}^k(d_i+1)-\sum_{i=2}^k(r_i+1)=
 \sum_{i=1}^kd_i-\sum_{i=2}^kr_i+1=n$, so the degree of the numerator $p(t)$ of the Hilbert series
 $\frac{p(t)}{(1-t)^n}$ of $k[F]$ is $n-\min\{ r_i\}-1$ so the projective dimension is $n-\min\{ r_i\}-2$,
 and the depth of $k[F]$ is $\min\{ r_i\}+2$ by the Auslander-Buchsbaum formula. We have
 $\dim k[F]=1+\max\{ d_i\}$, depth $k[F]=\min\{ r_i\}+2$, and $d_i>r_i$ for all $i$. The only possibility
 for $\dim k[F]={\rm depth}\, k[F]$ is that there is a $d$ such that $d_i=r_i+1=d$ for all $i$.

\bigskip\noindent
{\em Remark} The characterization of CM fat trees is not new. With another (more complicated)
proof it is found in \cite{Fr}.

\bigskip
In the remaining part of this section, we give some examples of results we can achieve.

\bigskip
Jacques has determined the Betti numbers of the complete bipartite graph $K_{m,n}$, \cite{Ja}.
The result is that the only nonzero Betti numbers are $\beta_{i,i+1}(K_{m,n})=
\sum_{j+l=i+1\atop j,l\ge1}{m\choose j}{n\choose l}$. We give an alternative proof.

\begin{theorem}
The edge ring of $K_{m,n}$ has a 2-linear resolution and $\beta_{i,i+1}(K_{m,n})=
{m+n\choose i+1}-{m\choose i+1}-{n\choose i+1}$.
\end{theorem}

\noindent
{\em Proof}  The edge ring of $K_{m,n}$ is the Stanley-Reisner ring of $K_m\sqcup K_n$, the
disjoint union of two complete graphs, so the resolution is 2-linear. The Hilbert series is
$\frac{1}{(1-t)^m}+\frac{1}{(1-t)^n}-1=\frac{(1-t)^n+(1-t)^m-(1-t)^{m+n}}{(1-t)^{m+n}}$.
Thus $\beta_{i,i+1}={n+m\choose i+1}-{n\choose i+1}-{m\choose i+1}$. 

\begin{cor}
$\sum_{j+l=i+1\atop j,l\ge1}{m\choose j}{n\choose l}={n+m\choose i+1}-{n\choose i+1}-{m\choose i+1}$.
\end{cor}

Also the complete multipartite graph $K_{n_1,\ldots,n_s}$ with parts of size $n_1,\ldots,n_s$ 
is treated in \cite{Ja}. The result there is
$$\beta_{i,i+1}(K_{n_1,\ldots,n_s})=
\sum_{l=2}^s(l-1)\sum_{\alpha_1+\cdots+\alpha_l=i+1\atop {\alpha_1,\dots,\alpha_l\ge1\atop j_1<\cdots<j_l}}
{n_{j_1}\choose\alpha_1}\cdots{n_{j_l}\choose\alpha_l}.$$
With our method we get

\begin{theorem}
$\beta_{i,i+1}(K_{n_1,\ldots,n_s})=\sum_{i=1}^s{N-n_i\choose i+1}-{N\choose i+1}$, where 
$N=\sum_{i=1}^sn_i.$
\end{theorem}

\noindent
{\em Proof}
The edge ring of $K_{n_1,\ldots,n_s}$ is the Stanley-Reisner ring of the disjoint union of
$K_{n_1},\ldots,K_{n_s}$ with Hilbert series $\sum_{i=1}^s\frac{1}{(1-t)^{n_i}}-(s-1)=
\frac{\sum_{i=1}^s(1-t)^{N-n_i}-(s-1)(1-t)^N}{(1-t)^N}$. 
The result follows as before.

\begin{cor}
$$\sum_{l=2}^s (l-1)\sum_{\alpha_1+\cdots+\alpha_l=i+1\atop{j_1<\cdots<j_l\atop\alpha_1,\dots,\alpha_l\ge1}}
{n_{j_1}\choose\alpha_1}\cdots{n_{j_l}\choose\alpha_l}=\sum_{i=1}^s{N-n_i\choose i+1}-{N\choose i+1}.$$
\end{cor}

In \cite{Ar-He-Hi} squarefree lexsegments ideals with $q$-linear resolution are studied. They show e.g.
that ideals generated by initial segment of squarefree monomials in degree $q$ in lexicographic order have linear resolution. 
For $q=2$ this means that ideals $L(a,b)$ generated by all squarefree monomials of degree 2 that are larger than or equal to $x_ax_b$
for some $(a,b)$ have 2-linear resolution.

\begin{theorem}
The Betti numbers of $k[x_1,\ldots,x_n]/L(a,b)$, $a\le b$, are $\beta_{i,i+1}=a{b\choose i+1}-a{b-1\choose i+1}-{a\choose i+1}$, $1\le i\le b-1$.
\end{theorem}

\noindent
{\em Proof}
$L(a,b)$ is the Stanley-Reisner ideal of a simplicial complex with maximal faces $\{1\},\{2\},\ldots,\{ a\},\{ a+1,a+2,\ldots,b\}$.
Thus the Hilbert series is $\frac{a}{(1-t)}+\frac{1}{(1-t)^{b-a}}-a=\frac{a(1-t)^{b-1}+(1-t)^a-a(1-t)^b}{(1-t)^b}$.

\begin{cor}
$a{b\choose i+1}-a{b-1\choose i+1}-{a\choose i+1}=\sum_{k=0}^{a-1}(k+1){k\choose i-1}+a\sum_{k=a}^{b-2}{k\choose i-1}$.
\end{cor}

\noindent
{\em Proof}
In \cite{Ar-He-Hi} it is shown that $\beta_{i,i+1}=\sum_{x_ix_j\ge x^ax^b}{j-2\choose i-1}$, $i\le j$.

\bigskip
Also final squarefree segment define ideals with linear resolution. For $d=2$ this means that ideals $F(a,b)=(\{ x_ix_j;x_ix_j\le x_ax_b\})$
have 2-linear resolution.

\begin{theorem}
The Hilbert series of $k[x_1,\ldots,x_n]/F(a,b)$, $a<b$, is 
$$\frac{n-b-1}{(1-t)^{N-a}}+\frac{b-a-1}{(1-t)^{N-a+1}}-\frac{b-a-2}{(1-t)^{N-a}}-\frac{n-b+1}{(1-t)^{N-a+1}},$$
where $N={n-a-1\choose 2}+n-b+1$.
\end{theorem}

\noindent{\em Proof}
$F(a,b)$ is the Stanley-Reisner ideal of a simplicial complex with maximal faces $\{1,2,\ldots,a-1,i\}$, $b\le i\le n$, and $\{1,2,\ldots,a,j\}$, $a+1\le j\le b-1$.

\bigskip\noindent
Now we treat some examples of Ferrer's ideals. A Ferrer's ideals can be indexed by a tableau $(\lambda_1,\lambda_2\ldots,\lambda_m)$,
$\lambda_1\ge\lambda_2\cdots\ge\lambda_m\ge1$. Here $$I_{(\lambda_1,\ldots,\lambda_m)}=
(x_1y_1,x_1y_2,\ldots,x_1y_{\lambda_1},x_2y_1,\ldots,x_2y_{\lambda_2},\ldots,x_my_1,\ldots,x_my_{\lambda_m}).$$
We start with an example from \cite{Co-Na}, namely $k[x_1,\ldots,x_5,y_1,\ldots,y_6]/I_{(6,4,4,2,1)}$. The fat tree with Stanley-Reisner ideal
$I_{(6,4,4,2,1)}$ can be constructed like this. (We denote the simplex on $a_1,\dots,a_k$ by $[a_1,\ldots,a_k]$.) 

Start with $[y_1,\ldots,y_6]$
and attach $[x_5,y_2,y_3,y_4,y_5,y_6]$ in $[y_2,y_3,y_4,y_5,y_6]$. Then attach\\
$[x_4,x_5,y_3,y_4,y_5,y_6]$ in $[x_5,y_3,y_4,y_5,y_6]$. Then
attach $[x_2,x_3,x_4,x_5,y_5,y_6]$ in $[x_4,x_5,y_5,y_6]$. Finally attach $[x_1,x_2,x_3,x_4,x_5]$ in $[x_2,x_3,x_4,x_5]$. 

The Hilbert series
of $k[x_1,\ldots,x_5,y_1,\ldots,y_6]/I_{(6,4,4,2,1)}$ is 
$$\frac{1}{(1-t)^6}+\frac{1}{(1-t)^6}-\frac{1}{(1-t)^5}+\frac{1}{(1-t)^6}-\frac{1}{(1-t)^5}
+\frac{1}{(1-t)^6}-\frac{1}{(1-t)^4}+\frac{1}{(1-t)^5}-\frac{1}{(1-t)^4}=$$
$$\frac{4(1-t)^5-(1-t)^6-2(1-t)^7}{(1-t)^{11}}.$$
Thus 
$$\beta_i(k[x_1,\ldots,x_5,y_1,\ldots,y_6]/I_{(6,4,4,2,1)})=2{7\choose i+1}+{ 6\choose i+1}-4{5\choose i+1}.$$ 

The formula
in \cite{Co-Na} gives 
$$\beta_i(k[x_1,\ldots,x_5,y_1,\ldots,y_6]/I_{(6,4,4,2,1)})={6\choose i}+{5\choose i}+{6\choose i}+{5\choose i}+{5\choose i}
-{5\choose i+1}=$$

$$2{6\choose i}+3{5\choose i}-{5\choose i+1}.$$
Denote a sequence $a,a,\ldots,a$ of length $b$ with $a^b$. A Ferrer's graph with tableau $\mu_1^{l_1},\ldots,\mu_k^{l_k}$ corresponds to
a simplicial complex with maximal faces $[x_{L-l_k-\cdots-l_i+1},\ldots,x_L,y_{\mu_i+1},\ldots,y_{\mu_1}]$,
where $L=\sum_{i=1}^kl_i$, for $i=0,\ldots,k$.
One could generalize the method above to any tableau, but the result becomes a bit complicated, so we refrain from doing so.
Instead we give some concrete examples and compare the formulas for the Betti numbers. 

\bigskip
First consider $I_{m^n}$ which is the edge ring of
the complete bipartite graph $K_{m,n}$.  Here we have yet another expression from \cite{Co-Na}, $\beta_{i,i+1}=\sum_{j=n}^{n+m-1}{j\choose i}-{m\choose i+1}$.
Thus we get 

\begin{cor}
${n+m\choose i+1}-{n\choose i+1}-{m\choose i+1}=\sum_{j=n}^{n+m-1}{j\choose i}-{m\choose i+1}$.
\end{cor}

\bigskip
Now consider $$I_ {n,1^{m-1}}=(x_1y_1,x_1y_2,\ldots,x_1y_n,x_2y_1,x_3y_1,\ldots,x_my_1).$$
The fat tree with this Stanley-Reisner ideal has a maximal face $[y_1,\ldots,y_n]$. Then another maximal face $[x_2,\ldots,x_m,y_2,\ldots,y_n]$ is attached in 
$[y_2,\ldots,y_n]$. Finally the maximal face $[x_1,\ldots,x_m]$ is attached in $[x_2,\ldots,x_m]$.

\begin{theorem}
The Hilbert series of $k[x_1,\ldots,x_m,y_1,\ldots,y_n]/I_{n,1^{m-1}}$ is 
$$\frac{1}{(1-t)^n}+\frac{1}{(1-t)^m}+\frac{1}{(1-t)^{m+n-2}}-\frac{1}{(1-t)^{n-1}}-\frac{1}{(1-t)^{m-1}},$$ 
so the Betti numbers are 
$$\beta_i={n+1\choose i+1}+{m+1\choose i+1}-{n\choose i+1}-{m\choose i+1}-{2\choose i+1}.$$
\end{theorem}

\begin{cor}
$${n+1\choose i+1}+{m+1\choose i+1}-{n\choose i+1}-{m\choose i+1}-{2\choose i+1}={n\choose i}+\sum_{j=2}^m{j\choose i}-{m+1\choose i+1}.$$
\end{cor}

\noindent
{\em Proof} The right hand side is the formula for $\beta_i$ from \cite{Co-Na}.

\bigskip
We can do the same in higher dimension. Consider the ideal 
$$I(m,n,o)=(x_1y_i,x_1z_j,y_1x_k,y_1z_j,z_1x_k,z_1y_i,2\le k\le m,2\le i\le n,2\le j\le o).$$
In this case the complex has four maximal faces,
$$[x_2,\ldots,x_m,y_2,\ldots,y_n,z_2,\ldots,z_o],[x_1,\ldots,x_m],[y_1,\ldots,y_n],[z_1,\ldots,z_o].$$
The three last are attached in $[x_2,\ldots,x_m]$, $[y_2,\ldots,y_n]$, $[z_2,\ldots,z_o]$, respectively. Then $k[x_1,\ldots,x_m,y_1,\ldots,y_n,z_1\ldots,z_o]/I(m,n,o)$
has a 2-linear resolution. The Hilbert series of $k[x_1,\ldots,x_m,y_1,\ldots,y_n,z_1\ldots,z_o]/I(m,n,o)$
is 
$$\frac{1}{(1-t)^{m+n+o-3}}+\frac{1}{(1-t)^m}+\frac{1}{(1-t)^n}\frac{1}{(1-t)^o}-\frac{1}{(1-t)^{m-1}}-\frac{1}{(1-t)^{n-1}}-\frac{1}{(1-t)^{o-1}},$$ 
so the Betti numbers are 
$$\beta_i={m+n+1\choose i+1}+{m+o+1\choose i+1}+{n+o+1\choose i+1}-{m+n\choose i+1}-{m+o\choose i+1}+{n+o\choose i+1}
-{3\choose i+1}.$$

\bigskip\noindent
{\em Example} Consider the Ferrer's ideal with tableau $n,n-1,n-2,\ldots,1$. Using the formula in  \cite{Co-Na} the Betti numbers are $n{n\choose i}-{n\choose i+1}$.
The fat tree giving this ideal is constructed like this. Start with $[y_1,\ldots,y_n]$ and attach $[x_n,y_2,\ldots,y_n]$ in $[y_2,\ldots,y_n]$. Then attach
$[x_{n-1},x_n,y_3,\ldots,y_n]$ in $[x_n,y_3,\ldots,y_n]$. Continue like this until: Attach $[x_2,\ldots,x_n,y_n]$ in $[x_3,\ldots,x_n,y_n]$, and finally
attach $[x_1,\ldots,x_n]$ in $[x_2,\ldots,x_n]$. The Hilbert series is $\frac{n+1}{(1-t)^n}-\frac{n}{(1-t)^{n-1}}$, so the Betti numbers are
$\beta_{i,i+1}=n{n+1\choose i+1}-(n+1){n\choose i+1}$ and we get\\ $n{n+1\choose i+1}-(n+1){n\choose i+1}=n{n\choose i}-{n\choose i+1}$.

\section{Uniform fat forests}
We will now concentrate on ``uniform" fat forests. If $F=F_1\cup\cdots\cup F_k$, where $F_i$
is a $d$-simplex for $i=1,\ldots,k$ and $(F_1\cup\cdots\cup F_{i-1})\cap F_i$ is an $r$-simplex for
$i=2,\ldots,k$, we call $F$ a $(d,r)$-forest.

\begin{cor}
The Hilbert series of a $(d,r)$-forest $\Delta$ with $k$ facets is 
$\frac{k}{(1-t)^{d+1}}-\frac{k-1}{(1-t)^{r+1}}$. The depth of $k[\Delta]$ is $r+2$, so $k[\Delta]$
is CM if and only if $r=d-1$. We have that the Betti numbers 
$b_i(k[\Delta])=b_{i,i+1}(k[\Delta])=|k{(k-1)(d-r)\choose i+1}-(k-1){k(d-r)\choose i+1}|$. 
\end{cor}

For a simplicial complex $\Delta$ on $[{\bf n}]$, the Alexander dual $\Delta^\vee$ is defined as 
$\{ F;F^c\notin\Delta\}$, where $F^c=[{\bf n}]\setminus F$. Alexander duals are well described in
\cite{Br-He}, \cite{Ba}, \cite{Ea-Re}, \and \cite{Te}. We will now determine the Betti
numbers of the Alexander dual to fat forests.

\begin{theorem}
Let $\Delta$ be a $(d,r)$-forest with $k$ facets.
The nonzero Beitti numbers of $k[\Delta^\vee]$ are 
$b_0=b_{0,0}=1$, $b_1=b_{1,(k-1)(d-r)}=k$ and $b_2=b_{2,k(d-r)}=k-1$. All $k[\Delta^\vee]$
are CM, and $k[\Delta^\vee]$ has a linear resolution if and only if $\Delta$ is a $(d,d-1)$-forest.
\end{theorem}

\noindent
{\em Proof} We have that $\prod_{j=1}^mx_{i_j}$ is a minimal generator of $I$, where
$k[\Delta^\vee]=k[x_1,\ldots,x_n]/I$, if and only if $[{\bf n}]\setminus\{ i_1,\ldots,i_m\}$
is a facet in $\Delta$. Thus $\Delta^\vee$ has $k$ minimal generators of degree
$k(d-r)+r+1-(d+1)=(k-1)(d-r)$, so $b_1(k[\Delta^\vee])=b_{1,(k-1)(d-r)}(k[\Delta^\vee])=k$.
The regularity of $k[\Delta]$ equals the projective dimension of $k[\Delta^\vee]$ and
vice versa, \cite{Te}, \cite{Te1}. The projective dimension of $k[\Delta]$ is $k(d-r)-1$ and
$b_{k(d-r)-1,k(d-r)}=k-1$. This is an extremal Betti number for $k[\Delta]$. (We have that
$b_{k,l}(R)$ is an extremal Betti number for $R$ if $b_{i,j}(R)=0$ for 
$\{(i,j);i\ge k,j\ge l\}\setminus\{k,l)\}$.) Then $b_{2,k(d-r)-1}(\Delta^\vee)=k-1$ is an
extremal Betti number, see \cite{Ba}. Since $\sum_{i=0}^2(-1)^ib_i=0$ we have $b_2=k-1$,
so $b_{2,j}=0$ if $j\ne k(d-r)-1$. We have that $k[\Delta^\vee]$ has a linear resolution if and
only if $k[\Delta]$ is CM, \cite{Ea-Re}.

\medskip
{\bf Note} No datasets were generated or analysed during the current study.

\end{document}